\documentclass{elsart}

\usepackage{amssymb,latexsym,amsfonts}

\newcommand{\tensor}{\otimes}

\renewcommand{\sec}{\Sigma_1}

\renewcommand{\c}{\operatorname{c}}

\newcommand{\picard}{\operatorname{Pic}}

\renewcommand{\H}{\operatorname{H}}
\newcommand{\h}{\operatorname{h}}

\newcommand{\wtp}{\widetilde{\P^n}}

\newcommand{\wts}{\widetilde{\Sigma}_1}

\newcommand{\ses}[3]{0\rightarrow#1\rightarrow#2
   \rightarrow#3\rightarrow0}

\newcommand{\F}{{\mathcal F}}
\newcommand{\G}{{\mathcal G}}

\newcommand{\I}{{\mathcal I}}

\renewcommand{\O}{{\mathcal O}}
\renewcommand{\P}{{\mathbb{P}}}

\newcommand{\N}{{\operatorname{N}}}

\newenvironment{proof}{\par \medskip \noindent
{\sc Proof:}}{}

\newtheorem{lemma}[thm]{Lemma}         
 
\renewenvironment{rem}[2]{\refstepcounter{thm} \label{#2} 
\par \medskip \noindent {\bf #1 \thethm .}}{\par \medskip}

\begin{document}

\begin{frontmatter}

\pagenumbering{arabic}

\title{Regularity and Normality of the Secant Variety to a Projective Curve}

\author{Peter Vermeire}

\ead{verme1pj@cmich.edu}

\address{Department of Mathematics, 214 Pearce, Central Michigan
University, Mount Pleasant MI 48859}



\date{\today}

\begin{abstract} 
For a smooth curve of genus $g$ embedded by a line bundle of degree at least $2g+3$ we show that the ideal sheaf of the secant variety is $5$-regular.  This bound is sharp with respect to both the degree of the embedding and the bound on the regularity.  Further, we show that the secant variety is projectively normal for the generic embedding of degree at least $2g+3$.

\end{abstract}

\end{frontmatter}


\section{The Result}

We work over an algebraically closed field of characteristic $0$. 
Recall the classical theorem of Castelnuovo:
\begin{thm}
Let $C\subset\P^n$ be a linearly normal embedding of a smooth curve of genus $g$ by a line bundle $L$ with $\c_1(L)\geq 2g+1$. Then $\I_C$ is $3$-regular and (equivalently) $C\subset\P^n$ is projectively normal.
{\nopagebreak \hfill $\Box$ \par \medskip}
\end{thm}

The following extension was proved for $a=2$ by J. Rathmann \cite{rathmann} and was proved in general by the author \cite[4.2]{vermeirevanishing} (see also \cite{bertramvanishing}).
\begin{thm}\label{extension}
Let $C\subset\P^n$ be a linearly normal embedding of a smooth curve of genus $g$ by a line bundle $L$ with $\c_1(L)\geq 2g+3$. Then $\I^a_C$ is $(2a+1)$-regular.
{\nopagebreak \hfill $\Box$ \par \medskip}
\end{thm}

Considering $C$ to be the zeroth secant variety to itself, and denoting the first secant variety by $\sec$, in this work we obtain what is perhaps a more natural extension.  
\begin{thm}\label{weakmain}
Let $C\subset\P^n$ be a linearly normal embedding of a smooth curve of genus $g$ by a line bundle $L$ with $\c_1(L)\geq 2g+3$. Then $\I_{\sec}$ is $5$-regular and at least for the generic such embedding $\sec\subset\P^n$ is projectively normal. 
{\nopagebreak \hfill $\Box$ \par \medskip}
\end{thm}

We expect that $\sec\subset\P^n$ is always projectively normal under these hypotheses.  The most significant difficulty in the proof of Theorem~\ref{weakmain} is that $\sec$ has rational singularities if and only if $C$ itself is rational.  Thus the usual technique \cite{bel} of blowing up a variety and studying line bundles on the blow-up rather than ideal sheaves on $\P^n$ requires significant care.

\begin{rem}{Remark}{sharp}
Note that in $\P^4$, the secant variety to a non-degenerate elliptic curve of degree $5$ is a quintic hypersurface, and that the secant variety to a non-degenerate genus $2$ curve of degree $6$ is an octic hypersurface; hence Theorem~\ref{weakmain} is sharp.
{\nopagebreak \hfill $\Box$ \par \medskip}
\end{rem}

Because the $k$th secant variety $\Sigma_k$ to an elliptic normal curve $C\subset\P^{2+2k}$ is a  hypersurface of degree $2k+3$, Theorem~\ref{weakmain} and Remark~\ref{sharp} suggest:
\begin{conj}\label{bestpossible}
Let $C\subset\P^n$ be a linearly normal embedding of a smooth curve of genus $g$ by a line bundle $L$.  If $\c_1(L)\geq 2g+1+2k$, $k\geq0$, then $\I_{\Sigma_k}$ is $(2k+3)$-regular and $\Sigma_k$ is projectively normal.
{\nopagebreak \hfill $\Box$ \par \medskip}
\end{conj}

We combine this with a previous conjecture of the author \cite[3.10]{vermeireflip2} to form the following Green-Lazarsfeld type conjecture.  Following \cite{EGHP} we say that a closed projective scheme $X\subset\P^n$ satisfies $\N_{k,p}$ if the ideal of $X$ is generated by forms of degree $k$ and the syzygies are linear for $p-1$ steps.
\begin{conj}\label{green-laz}
Let $C\subset\P^n$ be a linearly normal embedding of a smooth curve of genus $g$ by a line bundle $L$.  If $\c_1(L)\geq 2g+1+p+2k$, $p,k\geq0$, then $\Sigma_k$ satisfies $\N_{k+2,p}$.
{\nopagebreak \hfill $\Box$ \par \medskip}
\end{conj}

For $k=0$ this is the famous result of \cite{mgreen} (see also \cite{laz}).  For $g=0$ this seems well-known.  For $g=1$ this was proven in \cite{vbh} and in \cite{fisher}.  For $g=2$, calculations done by J. Sidman \cite{sidman} with Macaulay 2 \cite{mac2} support the conjecture for $\c_1(L)\leq 13$.

\section{The Proof}
We denote the $i$th secant variety to an embedded projective curve $C\subset\P^n$ by $\Sigma_i$.  Note that $\Sigma_0=C$.

A line bundle $L$ on a curve $C$ is said to be {\em $k$-very ample}
if $h^0(C,L(-Z))=h^0(C,L)-k$ for all $Z\in S^kC$.
We recall (the first stages of) Bertram's `Terracini Recursiveness' construction, which provides the geometric framework for our results. 

\begin{thm}\cite[Theorem 1]{bertram1}\label{terracini}
Let $C\subset B_0=\P(H^0(C,L))$ be a smooth curve embedded by a line
bundle $L$.
Suppose $L$ is $4$-very ample and consider the birational morphism
$f:B_2 \rightarrow B_0$
which is a composition of the following blow-ups:

$f^1:B_1\rightarrow B_0$ is the blow up of $B_0$ along $\Sigma_0$

$f^2:B_2 \rightarrow B_1$ is the blow up along the proper transform of $\Sigma_1$

Then, the proper transform of $\Sigma_1$ in $B_1$ is smooth and irreducible, transverse to the exceptional divisor, so in particular $B_2$ is smooth.  Let $E_i$ be the proper transform in $B_i$ of each $f^i$-exceptional divisor. 

(Terracini recursiveness) 
Suppose $x\in \Sigma_1\setminus C$. Then the fiber $f^{-1}(x) \subset B_2$ is naturally isomorphic to $\P(H^0(C,L(-2Z)))$, where $Z$ is the unique divisor of degree $2$ whose span contains $x$. If $x\in C$ the fiber $f^{-1}(x) \subset E_1 \subset B_2$ is isomorphic to the blow up of $\P(H^0(C,L(-2x)))$ along the image of $C$ embedded by $L(-2x)$.
{\nopagebreak \hfill $\Box$ \par \medskip}
\end{thm}

\begin{lemma}\label{tech}\cite[3.2]{vermeirevanishing}
Hypotheses and notation as above:
\begin{enumerate}
\item $\Sigma_1\subset B_0$ is normal.
\item $f_*\O_{B_2}=\O_{B_0}$ and $R^jf_*\O_{B_2}=0$ for $j\geq 1$.
\item $f_*\O_{E_i}=\O_{\Sigma_{i-1}}$ for $i=1,2$.
{\nopagebreak \hfill $\Box$ \par \medskip}
\end{enumerate}
\end{lemma}

Our proof proceeds by the well-known technique (Cf. \cite{bel}) of obtaining vanishings on the blow-ups, and then deducing vanishing statements on $B_0$.  In the case of a smooth variety $X\subset\P^n$, one has $\H^i(\wtp,\O(kH-aE))=\H^i(\P^n,\I_X^a(k))$.  The significant difficulty in the case of secant varieties is the following:
\begin{prop}\label{ratsing}
Let $C\subset\P^n$ be a $4$-very ample embedding of a smooth curve.  Then  $\sec$ has rational singularities if and only if $C$ is rational.
\end{prop}

\begin{proof}
Consider the blow up $f^1:B_1\rightarrow\P^n$.  By \cite[3.8]{vermeireflips1} we have $E_1\cap\wts=C\times C$ with the restriction $f^1:E_1\cap\wts\rightarrow C$ just projection onto one factor.  From the sequence
$$\ses{\I_{C\times C/B_1}}{\O_{B_1}}{\O_{C\times C}}$$
we see $R^2f^1_*\I_{C\times C/B_1}=\H^1(C,\O_C)\tensor\O_C$, which in turn implies $R^2f^1_*\I_{\wts/B_1}=\H^1(C,\O_C)\tensor\O_C$.  Because $R^2f^1_*\I_{\wts/B_1}=R^1f^1_*\O_{\wts}$, this completes the proof.
{\nopagebreak \hfill $\Box$ \par \medskip}
\end{proof}

Thus, in the notation of Theorem~\ref{terracini}, $R^1f_*\O_{E_2}\neq0$ for a non-rational curve and so the transfer of vanishing from $B_2$ to $B_0$ requires some care. In particular, $\H^2(B_1,\I_{\wts}(k))\neq0$ for all $k>>0$.  Our main result is:

\begin{thm}\label{reg}
Let $C\subset\P^n$, $n\geq5$, be a smooth curve embedded by a non-special line bundle $L$.  Assume that $C$ is linearly and cubically normal, and that there is a point $p\in C$ such that $L(2p)$ is $6$-very ample and $L(2p-2q)$ satisfies $K_2$ for all $q\in C$.  Then $\I_{\sec}$ is $5$-regular.
\end{thm}

Recall \cite{vermeireflips1} that an embedding $C\subset\P^n$ satisfies $K_2$ if $\I_C(2)$ is globally generated and the Koszul syzygies are generated by linear syzygies.

\begin{cor}\label{deg}
Let $C\subset\P^n$ be a smooth curve embedded by a line bundle of degree at least $2g+3$.  Then $\sec$ is $5$-regular.
\end{cor}

\begin{proof}
It is well-known \cite{mgreen} that a line bundle of degree at least $2g+3$ satisfies the hypotheses of Theorem~\ref{reg}.  We need only mention that a rational normal curve of degree $3$ has $\sec=\P^3$, that the secant variety to a rational normal curve of degree $4$ is a cubic hypersurface, and that the secant variety to an elliptic normal curve of degree $5$ is a quintic hypersurface.
{\nopagebreak \hfill $\Box$ \par \medskip}
\end{proof}

\begin{proof}(of Theorem~\ref{reg})
By Proposition~\ref{ratsing} we know that $f_*\O_{B_2}(-E_2)=\I_{\sec}$, that $R^2f_*\O_{B_2}(-E_2)=\H^1(C,\O_C)\tensor\O_C$, and that $R^if_*\O_{B_2}(-E_2)=0$ for all other $i$.

\underline{{\boldmath $\H^1(\P^n,\I_{\sec}(4))=0$}}: Let $\F$ be a coherent sheaf on $\P^n$, $\G$ a coherent sheaf on $B_2$.  From the 5-term sequence $$0\rightarrow \H^1(\P^n,\F\tensor f_*\G)\rightarrow \H^1(B_2,f^*\F\tensor\G)\rightarrow \H^0(\P^n,\F\tensor R^1f_*\G)\rightarrow \cdots$$
associated to the Leray-Serre spectral sequence we see that it is enough to show that $\H^1(B_2,\O_{B_2}(4H-E_2))=0$.  

Consider the sequence
$$\ses{\O_{B_2}(4H-E_1-E_2)}{\O_{B_2}(4H-E_2)}{\O_{E_1}(4H-E_2)}$$
Along the fibers $F_c$ of $E_1\rightarrow C$ we have $\O_{B_2}(4H-E_2)\tensor\O_{F_c}=\O_{F_c}(-E)$.  Therefore $R^if_*\O_{E_1}(4H-E_2)=0$ for $i\neq2$, hence $\H^i(E_1,\O_{E_1}(4H-E_2))=0$ for $i=0,1$, and so it suffices to show $\H^1(B_2,\O_{B_2}(4H-E_1-E_2))=0$.  By \cite[3.3]{vermeirevanishing} we know that $\H^1(B_2,\O_{B_2}(3H-E_1-E_2))=0$; pulling the Euler sequence on $\P^n$ up to $B_2$ we have
$$\ses{f^*\Omega^1_{\P^n}\tensor\O_{B_2}(4H-E_1-E_2)}{\bigoplus_{1}^{n+1}\O_{B_2}(3H-E_1-E_2)}{\O_{B_2}(4H-E_1-E_2)}$$
and it suffices to show $\H^2(B_2,f^*\Omega^1_{\P^n}\tensor\O_{B_2}(4H-E_1-E_2))=0$.

By non-specialty of $L$ together with cubic normality of the embedding we have $\H^1(B_2,\O_{B_2}(4H-E_1))=\H^2(B_2,\O_{B_2}(3H-E_1))=0$, hence $\H^2(B_2,f^*\Omega^1_{\P^n}\tensor\O_{B_2}(4H-E_1))=0$.  Finally, we show $\H^1(E_2,f^*\Omega^1_{\P^n}\tensor\O_{E_2}(4H-E_1))=\H^1(\wts,f^*\Omega^1_{\P^n}\tensor\O_{\wts}(4H-E_1))=0$.

From the sequence $$\ses{\O_{B_2}(3H-E_1-E_2)}{\O_{B_2}(3H-E_1)}{\O_{E_2}(3H-E_1)}$$
because we know that $\H^1(B_2,\O_{E_2}(3H-E_1))=0$ by cubic normality and that $\H^2(B_2,\O_{B_2}(3H-E_1-E_2))=0$ by \cite[3.3]{vermeirevanishing}, we have $\H^1(E_2,\O_{E_2}(3H-E_1))=0$.  Again working with the Euler sequence
$$\ses{f^*\Omega^1_{\P^n}\tensor\O_{\wts}(4H-E_1)}{\H^0(\P^n,\O(1))\tensor\O_{\wts}(3H-E_1)}{\O_{\wts}(4H-E_1)}$$
we show $$\H^0(\P^n,\O(1))\tensor\H^0(\wts,\O_{\wts}(3H-E_1))\rightarrow\H^0(\wts,\O_{\wts}(4H-E_1))$$ is surjective.

We recall \cite{bertram1},\cite{vermeireflips1} that $\O_{B_1}(2H-E_1)$ is globally generated and that the restriction of the induced morphism $\varphi$ to $\wts$ is a $\P^1$ bundle $\varphi:\wts\rightarrow S^2C$; furthermore, $\wts=\P_{S^2C}(\varphi_*\O_{\wts}(H))$.  In particular, there is a very ample line bundle $\O_{S^2C}(1)$ such that $\varphi^*\O_{S^2C}(1)=\O_{\wts}(2H-E_1)$.  It is easy to check, again with the Euler sequence, that $\varphi_*\H^0(\P^n,\O(1))\tensor\O_{\wts}\rightarrow E=\varphi_*\O_{B_1}(H)$ is surjective, hence $E$ is globally generated.  Applying $\varphi_*$ to the Euler sequence from the previous paragraph yields
$$\cdots\rightarrow\varphi_*\H^0(\P^n,\O(1))\tensor E\tensor\O_{S^2C}(1)\rightarrow S^2E\tensor\O_{S^2C}(1)\rightarrow0$$
By global generation of $E$ and of $\O_{S^2C}(1)$, and the vanishing of all higher direct images, we have $$\H^0(\P^n,\O(1))\tensor\H^0(\wts,\O_{\wts}(3H-E_1))\rightarrow\H^0(\wts,\O_{\wts}(4H-E_1))$$ is surjective.

\underline{{\boldmath $\H^2(\P^n,\I_{\sec}(3))=0$}}:  Again by \cite[3.3]{vermeirevanishing} we know that $\H^i(B_2,3H-E_1-E_2)=0$ for $i\geq1$.  Therefore we have
\begin{eqnarray*}
\H^2(\O_{B_2}(3H-E_2))&=&\H^2(\O_{E_1}(3H-E_2))\\
&=&\H^0(C,R^2f_*\O_{E_1}(3H-E_2))\\
&=&\H^0(C,R^2f_*\O_{B_2}(3H-E_2))
\end{eqnarray*}
It is straightforward to check that $E_2^{2,0}=E_{\infty}^{2,0}$, therefore because the edge map $\H^2(B_2,\O_{B_2}(3H-E_2))\rightarrow \H^0(C,R^2f_*\O_{B_2}(3H-E_2))$ is a quotient \cite[5.2.6]{weibel}, this implies
$\H^2(\P^n,f_*\O_{B_2}(3H-E_2))=\H^2(\P^n,\I_{\sec}(3))=0$.

\underline{{\boldmath $\H^3(\P^n,\I_{\sec}(2))=0$}}:  The fact that $\H^2(\wts,\O_{\wts}(2H-E_1))=0$ is contained in the proof of \cite[3.6]{vermeirevanishing}; therefore $\H^3(B_2,\O_{B_2}(2H-E_1-E_2))=0$.  Consider
$$\ses{\O_{B_2}(2H-E_1-E_2)}{\O_{B_2}(2H-E_2)}{\O_{E_1}(2H-E_2)}$$
We have $R^if_*\O_{E_1}(2H-E_2)=0$ for $i\neq2$ and $R^2f_*\O_{E_1}(2H-E_2)=\H^1(C,\O_C)\tensor\O_C(2)$.  Thus $\H^3(E_1,\O_{E_1}(2H-E_2))=\H^1(C,R^2f_*\O_{E_1}(2H-E_2))=0$, and so $\H^3(B_2,2H-E_2)=0$.  Further, we also have $$\H^2(B_2,\O_{B_2}(2H-E_2))\rightarrow \H^2(E_1,\O_{E_1}(2H-E_2))$$
is surjective, but as above $$\H^2(E_1,\O_{E_1}(2H-E_2))=\H^0(C,R^2f_*\O_{B_2}(2H-E_2))$$
and, therefore, $E_2^{0,2}=E_{\infty}^{0,2}$.  This implies that $d_3:E_3^{0,2}\rightarrow E_3^{3,0}=E_2^{3,0}$ is the zero map, and so $E_2^{3,0}=E_{\infty}^{3,0}$.  Thus the vanishing $\H^3(B_2,2H-E_2)=0$ implies $\H^3(\P^n,\I_{\sec}(2))=0$.

\underline{{\boldmath $\H^4(\P^n,\I_{\sec}(1))=0$}}:
Finally, $\H^3(\wts,\O_{B_1}(H))=\H^3(S^2C,\varphi_*\O_{B_1}(H))=0$.  Therefore, as it is not hard to see $E_2^{4,0}=E_{\infty}^{4,0}$, we have $\H^4(B_2,H-E_2)=\H^4(\P^n,\I_{\sec}(1))=0$.
{\nopagebreak \hfill $\Box$ \par \medskip}
\end{proof}

From the first part of the proof, we obtain the following general statement:
\begin{prop}\label{linnormal}
Let $C\subset\P^n$ be a $4$-very ample embedding of a smooth curve.  Then $\h^1(\P^n,\I^2_C(k))\geq\h^1(\P^n,\I_{\sec}(k))$ for $k\leq 3$. In particular, if $C$ is also linearly normal then $\sec$ is linearly normal.
\end{prop}

\begin{proof}
The proof shows 
\begin{eqnarray*}
\h^1(\P^n,\I_{\sec}(k))&=&\h^1(B_2,\O_{B_2}(kH-E_2))\\
&=&\h^1(B_2,\O_{B_2}(kH-2E_1-E_2))\\
&\leq&\h^1(B_2,\O_{B_2}(kH-2E_1))
\end{eqnarray*}
where the vanishing $\H^0(E_2,\O_{B_2}(kH-2E_1))=0$ for $k\leq 3$ comes from the structure of $\wts$ as a $\P^1$-bundle over $S^2C$.  The fact that $\h^1(\P^n,\I^2_C(1))=0$ for $C$ linearly normal is \cite[1.3.2]{wahl}.   
{\nopagebreak \hfill $\Box$ \par \medskip}
\end{proof}

\begin{rem}{Remark}{higherdim}
By the extension of Bertram's ``Terracini Recursiveness'' to higher dimensions \cite{vermeireflips1}, Proposition~\ref{linnormal} also holds for a $4$-very ample embedding of any projective variety $X\subset\P^n$ as long as the embedding satisfies $K_2$.
{\nopagebreak \hfill $\Box$ \par \medskip}
\end{rem}

It is worth pointing out that we do not always get equality in Proposition~\ref{linnormal}:
\begin{rem}{Example}{notfork=2}
Let $C\subset\P^4$ be a rational normal curve.  Then $\sec$ is a cubic hypersurface, hence $\H^1(\P^4,\I_{\sec}(2))=0$.  However, we can compute directly
that $\h^1(\P^n,\I^2_C(2))=3$.
{\nopagebreak \hfill $\Box$ \par \medskip}
\end{rem}

\begin{thm}
Let $C\subset\P^n$ be a smooth curve embedded by a non-special line bundle $L$.  Assume that $C$ is linearly and cubically normal, and that there is a point $p\in C$ such that $L(2p)$ is $6$-very ample and $L(2p-2q)$ satisfies $K_2$ for all $q\in C$.  Then for the general $q\in C$, the secant variety to $C$ under the embedding by $L(2p-2q)$ is projectively normal.
\end{thm}

\begin{proof}
Under these hypotheses, $\h^1(\P^n,\I^2_C(3))=0$ by \cite[3.10]{vermeirevanishing}, therefore we have $\H^1(\P^n,\I_{\sec}(3))=0$ by Proposition~\ref{linnormal}.
By Theorem~\ref{reg} and Proposition~\ref{linnormal}, we are left to show $\H^1(\P^n,\I_{\sec}(2))=0$.  This follows in general from \cite[3.5,3.9]{vermeirevanishing}.  In particular, the hypotheses allow us to construct a sequence of blow-ups $$f:B_3\rightarrow B_2\rightarrow B_1\rightarrow\P\Gamma(C,L)$$ where $B_3\rightarrow B_2$ is the blow up of the proper transform of $\Sigma_2$.  It is shown there that $R^1f_*\O_{E_1}(kH-2E_1-2E_2-E_3)=0$, which implies the generic vanishing of $\h^1(B_2,\O_{B_2}(2H-2E_2-E_2)$.  In order to get vanishing for ALL $q$ using this technique, one would need to show  $R^if_*\O_{E_3}(kH-2E_1-2E_2-E_3)=0$ for $i\geq2$.
{\nopagebreak \hfill $\Box$ \par \medskip}
\end{proof}

\begin{cor}
Let $C$ be a smooth curve.  For the generic $L\in\picard^kC$, $k\geq2g+3$, the secant variety $\sec\subset\P\Gamma(C,L)$ is projectively normal.
{\nopagebreak \hfill $\Box$ \par \medskip}
\end{cor}

\end{document}